%% file: 10-28.tex
\documentclass[a4,10pt]{article}
\usepackage{amsmath,amscd,amssymb}
\usepackage[margin=.8in]{geometry}

\input{notation}

\input{new_theorem}

\begin{document}

\title{Semi-flat metrics of the moduli spaces of Higgs bundles\\
in the non-zero degree case}

\author{Takuro Mochizuki\thanks{Research Institute for Mathematical Sciences, Kyoto University, Kyoto 606-8502, Japan, takuro@kurims.kyoto-u.ac.jp}}
\date{}
\maketitle

\begin{abstract}
We study horizontal deformations of 
a Higgs bundle whose spectral curve is smooth.
It allows us to define a natural integrable connection
of the Hitchin fibration on the locus where the spectral curves are smooth.
Then, in the non-zero degree case,
we introduce the semi-flat metric,
and compare the asymptotic behaviour of the semi-flat metric
and the Hitchin metric
along the ray $(E,t\theta)$ $(t\to\infty)$.

\vspace{.1in}
 \end{abstract}

\section{Introduction}

\subsection{Moduli space of stable Higgs bundles}

Let $X$ be a compact Riemann surface.
Let $\nbigm_{H}(X,n,d)$ denote the moduli space of
stable Higgs bundles of rank $n$ and degree $d$ on $X$.
By the theorem of Hitchin and Simpson
\cite{Hitchin-self-duality, s1},
it has a naturally defined hyperk\"ahler metric
$g_H^{X,n,d}$, called the Hitchin metric.
(We shall recall the construction in
\S\ref{subsection;25.1.22.10}.)

Let $\nbigm'_H(X,n,d)\subset\nbigm_H(X,n,d)$ 
denote the Zariski open subset
of Higgs bundles whose spectral curves are smooth.
We set $\vecA_{X,n}:=\bigoplus_{j=1}^nH^0(X,K_X^{j})$.
Let $\Phi_{n,d}:\nbigm_H(X,n,d)\to\vecA_{X,n}$
denote the Hitchin fibration (see \S\ref{subsection;25.1.22.12}).
We set $\vecA_{X,n}'=\Phi_{n,d}(\nbigm_H(X,n,d)')$.
The restriction $\Phi_{n,d}:\nbigm'_H(X,n,d)\to \vecA'_{X,n}$
is locally a principal torus bundle (see \S\ref{subsection;24.10.30.11}).

\subsection{Semi-flat metric in the case $d=0$}
\label{subsection;25.1.22.21}

In the case $d=0$,
$\nbigm'_H(X,n,0)$ has another hyperk\"ahler metric
$g_{\semiflat}^{X,n,0}$,
called the semi-flat metric.
Let us recall the construction.
There exists the Hitchin section
$\vecA_{X,n}\to\nbigm_H(X,n,0)$
of $\Phi_{n,0}$.
It induces a section
$\vecA'_{X,n}\to\nbigm'_H(X,n,0)$.
We obtain the integrable connection of
the locally principal torus bundle $\nbigm'_H(X,n,0)\to \vecA'_{X,n}$
for which the Hitchin section is horizontal.
For each $(E,\theta)\in\nbigm'_H(X,n,0)$,
we obtain the decomposition of the tangent space
into the vertical direction and the horizontal direction:
\begin{equation}
\label{eq;25.1.22.22}
 T_{(E,\theta)}\nbigm'_H(X,n,0)
=
 T_{(E,\theta)}\nbigm'_H(X,n,0)^{\ver}
 \oplus
 T_{(E,\theta)}\nbigm'_H(X,n,0)^{\hor}.
\end{equation}
Here,
$T_{(E,\theta)}\nbigm'_H(X,n,0)^{\ver}$
denotes the tangent space of the fiber
$T_{(E,\theta)}\Phi_{n,0}^{-1}(\Phi_{n,0}(E,\theta))$,
and 
$T_{(E,\theta)}\nbigm'_H(X,n,0)^{\hor}$
denotes the subspace determined by the integrable connection.
According to \cite[Corollary 3.33]{Mochizuki-Asymptotic-Hitchin-metric},
the horizontal space
$T_{(E,\theta)}\nbigm'_H(X,n,0)^{\hor}$ is Lagrangian.
Let $\Sigma_{\theta}$ denote the spectral curve of
the Higgs bundle $(E,\theta)$.
Because
$\Phi_{n,0}^{-1}(\Phi_{n,0}(E,\theta))$
is a connected component of the moduli space of holomorphic line bundles
on $\Sigma_{\theta}$,
there exists the following natural isomorphism
\[
 T_{(E,\theta)}\nbigm'_H(X,n,0)^{\ver}
 \simeq
 H^1(\Sigma_{\theta},\nbigo_{\Sigma_{\theta}}).
\]
As well known in \cite{Hitchin-self-duality},
each fiber of the Hitchin fibration $\Phi_{n,0}$ is Lagrangian,
and hence
\[
 T_{(E,\theta)}\nbigm'_H(X,n,0)^{\hor}
 \simeq
 \Hom\Bigl(
 T_{(E,\theta)}\nbigm'_H(X,n,0)^{\ver},
 \cnum
 \Bigr)
 \simeq
 H^0(\Sigma_{\theta},K_{\Sigma_{\theta}}).
\]
Therefore, there exists the natural isomorphism between
$T_{(E,\theta)}\nbigm_H'(X,n,0)$
and the space of harmonic $1$-forms on $\Sigma_{\theta}$.
The natural $L^2$-metric on the $1$-forms on $\Sigma_{\theta}$
induces
a metric $g^{X,n,0}_{\semiflat|(E,\theta)}$
on $T_{(E,\theta)}\nbigm_H'(X,n,0)$.
In this way, we obtain
a metric $g_{\semiflat}^{X,n,0}$ on $\nbigm'_H(X,n,0)$.
According to \cite{Freed}, it is hyperk\"ahler.
(See also \cite[Proposition 4.6]{Mochizuki-Asymptotic-Hitchin-metric}
and Proposition \ref{prop;25.1.22.100}.)

We can think that $g_{\semiflat}^{X,n,0}$
is easier compared with $g_H^{X,n,0}$
in the sense that $g_{\semiflat}^{X,n,0}$ is determined by
the geometry of $\nbigm_H'(X,n,0)$.

\subsection{Comparison of the Hitchin metric
and the semi-flat metric in the case $d=0$}

In \cite{Gaiotto-Moore-Neitzke},
Gaiotto, Moore and Neitzke proposed a beautiful conjectural description
of the hyperk\"ahler structure of $\nbigm_H(X,n,d)$.
In particular, it implies the following estimate
for any $(E,\theta)\in\nbigm_H(X,n,0)$ with respect to
$(g_{\semiflat}^{X,n,0})_{|(E,t\theta)}$ as $t\to\infty$:
\begin{equation}
\label{eq;25.1.22.20}
 \bigl(
 g_H^{X,n,0}-g^{X,n,0}_{\semiflat}
 \bigr)_{|(E,t\theta)}
 =O\bigl(\exp(-\epsilon t)\bigr).
\end{equation}
Here, $\epsilon$ denotes a positive number depending on
$(E,\theta)$.
The estimate (\ref{eq;25.1.22.20}) was studied in
\cite{Fredrickson2, Dumas-Neitzke, MSWW2}
and later in
\cite{Mochizuki-Asymptotic-Hitchin-metric,Mochizuki-improve}.

\subsection{The main result of this paper}

It is natural to study an analogue of (\ref{eq;25.1.22.20})
in the case $d\neq 0$.
In the case $d=0$,
as explained in \S\ref{subsection;25.1.22.21},
we obtain the integrable connection
of $\nbigm'_{H}(X,n,d)\to\vecA'_{X,n}$
by using the Hitchin section.
It induces the decomposition of the tangent space
(\ref{eq;25.1.22.22}),
which we used to construct $g_{\semiflat}^{X,n,0}$.
It is not directly applied to the case $d\neq 0$.

We first study
the notion of horizontal deformations
(see \S\ref{subsection;25.1.22.32}).
It allows us to obtain a natural integrable connection of
the locally principal torus bundle
$\nbigm'_H(X,n,d)\to \vecA'_{X,n}$,
whose horizontal spaces are Lagrangian,
and the semi-flat metric $g_{\semiflat}^{X,n,d}$
as in the case $d=0$.
Then, we obtain an analogue of
(\ref{eq;25.1.22.20})
in the case $d\neq 0$.
\begin{thm}[Theorem
\ref{thm;25.1.22.40}]
For any $(E,\theta)\in \nbigm'_{H}(X,n,d)$,
there exists $\epsilon>0$
such that the following estimate holds
as $t\to\infty$  with respect to $(g_{\semiflat}^{X,n,d})_{|(E,t\theta)}$:
\begin{equation}
\label{eq;25.1.22.50}
 \bigl(
 g_H^{X,n,d}-g^{X,n,d}_{\semiflat}
 \bigr)_{|(E,t\theta)}
 =O\bigl(\exp(-\epsilon t)\bigr).
\end{equation}
\end{thm}

We shall explain how to reduce the estimate (\ref{eq;25.1.22.50})
to the estimate (\ref{eq;25.1.22.20}) in the case $d=0$
in \cite{Mochizuki-Asymptotic-Hitchin-metric}.

\paragraph{Acknowledgement}

I prepared this manuscript for the proceedings volume
on ``Advances in Complex Lagrangians,
Integrable systems and Quantization''.
I thank the organizers for the invitation.

I gave a talk about a part of the content of this manuscript
in the conference to celebrate 60th birthday of Michihiko Fujii.
I thank the organizers, Tatsuro Shimizu and Shinpei Baba for the opportunity
of the talk.
I thank Bin Xu for discussions and kindness during my stay in
Chinese University of Science and Technology.
I reminded this issue when I prepared the lectures in Hefei.

I am partially supported by
the Grant-in-Aid for Scientific Research (A) (No. 21H04429),
the Grant-in-Aid for Scientific Research (A) (No. 22H00094),
the Grant-in-Aid for Scientific Research (A) (No. 23H00083),
and the Grant-in-Aid for Scientific Research (C) (No. 20K03609),
Japan Society for the Promotion of Science.
I am also partially supported by the Research Institute for Mathematical
Sciences, an International Joint Usage/Research Center located in Kyoto
University.

\section{Horizontal deformations of line bundles in a general setting}

\subsection{Preliminaries}

Let $B$ be a multi-disc.
Let $f:Z\to B$ be a proper morphism of complex manifolds
such that each fiber is a compact Riemann surface.
In particular, each fiber is assumed to be connected.
For each $b\in B$,
we set $Z_b:=f^{-1}(b)$,
and let $\iota_b:Z_b\to Z$ denote the inclusion.
Because $Z$ is diffeomorphic to $B\times Z_{b}$
for any $b\in B$,
there exists a natural isomorphism
$H_1(Z,\seisuu)\simeq H_1(Z_{b},\seisuu)$.
For any coherent $\nbigo_{Z}$-module
$\nbigf$ flat over $f^{-1}(\nbigo_B)$,
we set $\nbigf_b:=\iota_b^{\ast}(\nbigf)$.
We also set
$\deg(\nbigf):=\deg(\iota_{b}^{\ast}\nbigf)$
which is independent of $b\in B$.

We set $U(1)=\bigl\{a\in\cnum\,\big|\,|a|=1\bigr\}$.
We identify $H^1(Z,U(1))$
with the group of homomorphisms
of abelian groups $H_1(Z,\seisuu)\to U(1)$.
For any $\rho\in H^1(Z,U(1))$,
we obtain the induced unitary flat line bundle
$\nbigl(\rho)$ on $Z$.
The underlying holomorphic line bundle is also denoted by
$\nbigl(\rho)$.

For any non-zero integer $e$,
we set $\mu_{e}=\bigl\{a\in\cnum\,|\,a^{e}=1\bigr\}\subset U(1)$.
The subgroup
$H^1(Z,\mu_e)\subset H^1(Z,U(1))$
is identified with the group of homomorphisms
of abelian groups $H_1(Z,\seisuu)\to \mu_e$.

The following lemma is well known.
\begin{lem}
Let $b\in B$.
For any line bundle $N$ on $Z_b$ of degree $0$,
there uniquely exists $\rho\in H^1(Z,U(1))$
such that $N\simeq \nbigl(\rho)_b$.
\hfill\qed
\end{lem}

\begin{lem}
\label{lem;24.10.29.3}
Let $L$ be any holomorphic line bundle on $Z$
with $\deg(L)=0$.
Then, we have either 
$f_{\ast}(L)=0$
or 
$f_{\ast}(L)\simeq\nbigo_B$.
\end{lem}
\pf
Because $B$ is a multi-disc,
it is enough to prove that
$f_{\ast}(L)$ is locally free sheaf of
rank $1$ or $0$.
We have only to check the claim locally around
any point of $B$.

Let $b\in B$.
Let $U_b$ be any small neighbourhood of $b$ in $B$
which is isomorphic to a multi-disc.
We obtain $Z_{U_b}:=Z\times_BU_b$,
and the induced morphism $f_{U_b}:Z_{U_b}\to U_b$.
The restriction of $L$ to $Z_{U_b}$
is denoted by $L_{U_b}$.
We may assume that there exists a holomorphic section
$h:U_b\to Z_{U_b}$ of $f_{U_b}$.
Let $H\subset Z_{U_b}$ denote the image of $h$.
By shrinking $U_b$,
we may assume that there exists $m\in\seisuu_{>0}$
such that $R^1(f_{U_b})_{\ast}\bigl(L_{U_b}(mH)\bigr)=0$.
Then, $(f_{U_b})_{\ast}(L_{U_b}(mH))$ is a locally free
$\nbigo_{U_b}$-module.
Because there exists a natural monomorphism
$(f_{U_b})_{\ast}(L_{U_b})
\to
(f_{U_b})_{\ast}\bigl(L_{U_b}(mH)\bigr)$,
we obtain that
$(f_{U_b})_{\ast}(L_{U_b})$
is a torsion-free $\nbigo_{U_b}$-module.

If $L_{c}\not\simeq\nbigo_{Z_{c}}$
for some $c\in U_b$,
there exists a neighbourhood $U_{c}$ in $U_b$
such that
$L_{c'}\not\simeq \nbigo_{Z_{c'}}$
for any $c'\in U_c$.
Because
$\dim H^0(Z_{c'},L_{c'})=0$
for $c'\in U_{c}$,
we obtain $f_{\ast}(L_{U_c})=0$.
Because $f_{\ast}(L_{U_b})$ is torsion-free,
we obtain 
$f_{\ast}(L_{U_b})=0$.
If $L_{c}\simeq \nbigo_{Z_c}$ for any $c\in U_b$,
we obtain that
$\dim H^0(Z_{c},L_{c})=1$ for any $c\in U_b$,
and that $\dim H^1(Z_{c},L_{c})$ are constant
for $c\in U_{b}$.
Hence, 
$f_{\ast}(L_{U_b})$
is a locally free $\nbigo_{U_b}$-module of rank $1$.
\hfill\qed

\subsection{Roots of line bundles}

We recall the following standard fact.

\begin{prop}
Let $M$ be any holomorphic line bundle on $Z$.
Let $e\in\seisuu\setminus\{0\}$ be any divisor of $\deg(M)$,
i.e.,
a non-zero integer such that $\deg(M)\in e\seisuu$.
\begin{itemize}
 \item There exists a holomorphic line bundle $M_1$
       with an isomorphism
       $M_1^{e}\simeq M$.
 \item Let $M_2$ be another holomorphic line bundle
       with an isomorphism
       $M_2^e\simeq M$.
       Then,
       there exists
       $\rho\in H^1(Z,\mu_e)$
       such that
       $M_2=M_1\otimes \nbigl(\rho)$.
 \item Let $b\in B$.
       For any holomorphic line bundle $N$ on $Z_{b}$
       with an isomorphism
       $N^{e}\simeq M_{b}$,
       there exists a holomorphic line bundle
       $N_B$ on $B$
       such that
       $(N_B)_{b}\simeq N$
       and
       $N_B^{e}\simeq M$.
\end{itemize} 
\end{prop}
\pf
We explain a proof just for the convenience of readers.
Let $\nbigo_Z^{\times}$ denote the sheaf of
nowhere vanishing holomorphic functions.
Recall that
the isomorphism classes of holomorphic line bundles
on $Z$ is classified by $H^1(Z,\nbigo_{Z}^{\times})$.
Let $\seisuu_{Z}$ denote the sheaf of
$\seisuu$-valued continuous functions on $Z$.
From the exact sequence of sheaves
$0\lrarr
\seisuu_{Z}
\lrarr
\nbigo_{Z}
\lrarr
\nbigo_{Z}^{\times}
\lrarr 0$,
we obtain the exact sequence of modules
\[
 H^1(Z,\nbigo_{Z})
 \stackrel{a_1}{\lrarr}
 H^1(Z,\nbigo_{Z}^{\times})
 \stackrel{a_2}{\lrarr}
 H^2(Z,\seisuu)
 \simeq \seisuu
 \lrarr
 H^2(Z,\nbigo_{Z}).
\]
Let $\alpha\in H^1(Z,\nbigo_Z^{\times})$ 
be the element corresponding to $M$.
If $a_2(\alpha)\neq 0$, which is mapped to $0$ in $H^2(Z,\nbigo_Z)$,
we obtain that $a_2$ is surjective 
because $H^2(Z,\nbigo_Z)$ is a $\cnum$-vector space.
Hence, 
there exists $\beta\in H^1(Z,\nbigo_{Z}^{\times})$
such that
$a_2(\alpha-e\beta)=0$.
If $a_2(\alpha)=0$,
by setting $\beta=\alpha$,
we obtain $a_2(\alpha-e\beta)=0$.
Because $H^1(Z,\nbigo_{Z})$ is
a $\cnum$-vector space,
there exists
$\gamma\in H^1(Z,\nbigo_{Z})$
such that
$a_1\bigl(e\cdot\gamma\bigr)=\alpha-e\cdot\beta$.
Let $M_1$ be a line bundle corresponding to
$\beta-a_1(\gamma)$.
Then, there exists an isomorphism
$M_1^{e}\simeq M$.

Let $M_2$ be another holomorphic line bundle on $Z$
with an isomorphism $M_2^{e}\simeq M$.
There exists the following exact sequence:
\[
 0\lrarr
  \mu_{e}
  \lrarr
  \nbigo_{Z}^{\times}
  \stackrel{b}{\lrarr}
  \nbigo_{Z}^{\times}
  \lrarr 0.
\]
Here,
the morphism $b$ is defined by $b(g)=g^e$.
We obtain the following induced exact sequence:
\[
 H^1(Z,\mu_{e})
 \stackrel{b_1}{\lrarr}
 H^1(Z,\nbigo_{Z}^{\times})
 \stackrel{b_2}{\lrarr}
 H^1(Z,\nbigo_{Z}^{\times})
 \stackrel{b_3}{\lrarr}
 H^2(Z,\mu_{e}).
\]
Let $\alpha_i\in H^1(Z,\nbigo_{Z}^{\times})$
be the elements corresponding to $M_i$.
Because
$M_1^{e}
\simeq
M_2^{e}$,
we obtain
$b_2(\alpha_1)=b_2(\alpha_2)$.
There exists
$\beta\in H^1(Z,\mu_{e})$
such that
$\alpha_1=b_1(\beta)\alpha_2$.
Let
$\rho:H_1(Z_{b},\seisuu)\to \mu_{e}$
be the homomorphism corresponding to $\beta$.
It is easy to see that
$M_1\simeq M_2\otimes \nbigl(\rho)$.

Let $b\in B$.
Let $N$ be a holomorphic line bundle on $Z_{b}$
with an isomorphism $N^{e}\simeq M_{b}$.
Because
$N^e\simeq (M_1)_{b}^e$,
there exists
$\rho\in H^1(Z,U(1))$
such that
$N\simeq (M_1)_{b}\otimes \nbigl(\rho)_{b}$.
We set $N_B=M_1\otimes\nbigl(\rho)$.
Then,
there exist desired isomorphisms
$(N_B)_{b}\simeq N$
and
$N_B^{e}\simeq M$.
\hfill\qed

\begin{rem}
Because $B$ is assumed to be multi-disc,
it is easy to observe $H^2(Z,\nbigo_Z)=0$ in the proof.
\hfill\qed
\end{rem}

\subsection{Horizontal line bundles with respect to
a prescribed holomorphic line bundle}
\label{subsection;24.10.29.2}

Let $M$ be a holomorphic line bundle on $Z$.
Assume that $\deg(M)\neq 0$.

\begin{df}
\label{df;24.10.29.1}
A holomorphic line bundle $L$ on $Z$
is called horizontal with respect to $M$
if the following holds. 
\begin{itemize}
 \item Let $M_1$ be a holomorphic line bundle on $Z$
       with an isomorphism
       $M_1^{\deg(M)}\simeq M$.
       Then,
       there exists
       $\rho\in H^1(Z,U(1))$
       such that
       $L\simeq M_1^{\deg(L)}\otimes\nbigl(\rho)$.
       This condition is independent of the choice of $M_1$.
       \hfill\qed
\end{itemize}
\end{df}

\begin{rem}
If $\deg(L)=0$,
the condition is independent of $M$.
\hfill\qed
\end{rem}

\begin{lem}
\label{lem;24.10.29.31}
Let $L_i$ $(i=1,2)$ be holomorphic line bundles
on $Z$.
If $L_i$ are horizontal with respect to $M$,
then $L_1^{n_1}\otimes L_2^{n_2}$ are also horizontal
for any $(n_1,n_2)\in \seisuu^2$.
\end{lem}
\pf
There exit $\rho_i\in H^1(Z,U(1))$
such that
$L_i\simeq M_1^{\deg(L_i)}\otimes \nbigl(\rho_i)$.
Because
\[
 L_1^{n_1}\otimes L_2^{n_2}
=M_1^{n_1\deg(L_1)+n_2\deg(L_2)}\otimes\nbigl(n_1\rho_1+n_2\rho_2),
\]
the claim is clear.
\hfill\qed

\begin{prop}
\label{prop;24.10.30.3}
Let $N$ be any holomorphic line bundle on $Z_b$
for some $b\in B$.
\begin{itemize}
 \item 
       There exists a holomorphic line bundle
       $N_B$ on $Z$
       with an isomorphism
       $(N_B)_{b}\simeq N$
       such that $N_B$ is horizontal with respect to $M$.
 \item Let $N_B'$ be another such horizontal line bundle
       on $Z$
       with an isomorphism
       $(N_B')_{b}\simeq N$.
       Then, there exists
       an isomorphism
       $N_B'\simeq N_B$
       whose restriction to
       $Z_{b}$
       equals
       the composition of
       $(N_B')_{b}\simeq N\simeq (N_B)_{b}$.
\end{itemize}
\end{prop}
\pf
Let $M_1$ be a holomorphic line bundle on $Z$
with an isomorphism $M_1^{\deg(M)}\simeq M$.
There exists
$\rho\in H^1(Z,U(1))$
such that
$N\simeq (M_1)^{\deg(N)}_{b}\otimes\nbigl(\rho)_{b}$.
We obtain a desired line bundle
$N_B=M_1^{\deg(N)}\otimes\nbigl(\rho)$
with an isomorphism
$(N_B)_{b}\simeq N$.

Let $N_B'$ 
be another horizontal holomorphic line bundle on $Z$
with an isomorphism 
$(N_B')_{b}\simeq N$.
There exists 
$\rho'\in H^1(Z,U(1))$
such that
$N_B'=M_1^{\deg(N)}\otimes\nbigl(\rho')$.
Because
$(N_B)_{b}\simeq (N_B')_{b}$,
we obtain
$\nbigl(\rho)_{b}\simeq\nbigl(\rho')_{b}$,
which implies $\rho=\rho'$.
Hence, there exists a desired isomorphism
$N_B'\simeq N_B$.
\hfill\qed

\begin{lem}
\label{lem;24.10.30.2}
Let $M_i$ $(i=1,2)$ be holomorphic line bundles on $Z$
such that
there exist non-zero integers $m_i$
such that $M_1^{m_1}\simeq M_2^{m_2}$.
Then, a holomorphic line bundle $N$ on $Z$
is horizontal with respect to $M_1$
if and only if it is horizontal with respect to $M_2$. 
\end{lem}
\pf
By definition,
it is easy to see that
$N$ is horizontal with respect to $M_i$
if and only if $N$ is horizontal with respect to $M_i^{m_i}$.
Then, the claim is clear.
\hfill\qed

\subsection{Pull back by a covering map}

Let $f':Z'\to B$ be
a proper morphism of complex manifolds
such that each fiber is a compact Riemann surface.
Suppose that there exists a holomorphic covering map
$\psi:Z'\to Z$ such that $f'=f\circ\psi$.
For any 
$\rho\in H^1(Z,U(1))$,
we obtain
$\psi^{\ast}(\rho)\in H^1(Z',U(1))$.
There exists a natural isomorphism
$\psi^{\ast}\nbigl(\rho)
 \simeq
\nbigl(\psi^{\ast}\rho)$.
Let $\nbigk$ denote 
the kernel of
$\psi^{\ast}:
H^1(Z,U(1))
\to
H^1(Z',U(1))$.

\begin{lem}
$\nbigk$ is a finite abelian group.
\end{lem}
\pf
Let $e$ denote the degree of the covering map $\psi$.
Let $I\subset H_1(Z,\seisuu)$
denote the image of
$H_1(Z',\seisuu)\to H_1(Z,\seisuu)$.
For any $\gamma\in H_1(Z,\seisuu)$,
there exists
a positive integer $m(\gamma)$
such that $0<m(\gamma)\leq e$
and $m(\gamma)\gamma\in I$.
Hence, there exists a positive integer $m$ such that
$m\cdot H_1(Z,\seisuu)\subset I$.
If $\rho\in\nbigk$,
the restriction of $\rho$ to $m\cdot H_1(Z,\seisuu)$
is trivial,
and hence we obtain
$\nbigk\subset H^1(Z,\mu_m)$.
\hfill\qed

\begin{lem}
Let $L$ be a holomorphic line bundle on $Z$
such that $\psi^{\ast}(L)\simeq\nbigo_{Z'}$.
Then, there exists $\rho\in\nbigk$
such that $L\simeq \nbigl(\rho)$. 
\end{lem}
\pf
For each $b\in B$,
there exists
$\rho(b)\in\nbigk$
such that
$L_b\simeq \nbigl(\rho(b))_b$.
Because the dependence of $\rho(b)$ on $b\in B$
is continuous,
we obtain that there exists $\rho\in\nbigk$
such that
$L_b\simeq \nbigl(\rho)_b$
for any $b\in B$.
We set
$\Ltilde=L\otimes\nbigl(\rho)^{-1}$.
By using Lemma \ref{lem;24.10.29.3},
we obtain
$f_{\ast}(\Ltilde)\simeq \nbigo_B$.
By taking a nowhere vanishing global section of
$f_{\ast}(\Ltilde)$,
we obtain an isomorphism $L\simeq \nbigl(\rho)$.
\hfill\qed

\begin{cor}
\label{cor;24.10.29.4}
Let $L$ be a line bundle on $Z$
such that $\psi^{\ast}(L)\simeq \nbigo_{Z'}$
and $L_{b}\simeq \nbigo_{Z_{b}}$ for some $b\in B$.
Then, $L\simeq\nbigo_{Z}$.
\hfill\qed 
\end{cor}

Let $M$ be a holomorphic line bundle on $Z$
such that $\deg(M)\neq 0$
as in \S\ref{subsection;24.10.29.2}.
\begin{lem}
\label{lem;24.10.29.10}
A holomorphic line bundle $L$ on $Z$
is horizontal with respect to $M$
if and only if
$\psi^{\ast}(L)$ is horizontal 
with respect to $\psi^{\ast}(M)$.
\end{lem}
\pf
Let $e\in\seisuu_{>0}$
denote the degree of the covering map $\psi$.
Note that for any line bundle $L$ on $Z$,
we have
$\deg(\psi^{\ast}(L))=e\deg(L)$.
Let $M_1$ be a holomorphic line bundle on $Z$
with an isomorphism $M_1^{\deg(M)}\simeq M$.
Let $M_2'$ be a holomorphic line bundle on $Z'$
with an isomorphism
$(M_2')^e\simeq \psi^{\ast}(M_1)$.

If $L$ is horizontal,
there exists
$\rho\in H^1(Z,U(1))$
such that
$L\simeq M_1^{\deg(L)}\otimes \nbigl(\rho)$.
We obtain
\[
\psi^{\ast}(L)
\simeq
\psi^{\ast}(M_1)^{\deg(L)}\otimes
\psi^{\ast}(\nbigl(\rho))
\simeq
(M_2')^{\deg(\psi^{\ast}(L))}
\otimes
 \nbigl(\psi^{\ast}(\rho)).
\]
Hence, $\psi^{\ast}(L)$ is horizontal.

Conversely,
suppose that $\psi^{\ast}(L)$ is horizontal.
There exists a homomorphism
$\rho'\in H^1(Z',U(1))$
such that
$\psi^{\ast}(L)\simeq
(M'_2)^{\deg(\psi^{\ast}(L))}\otimes \nbigl(\rho')$,
i.e.,
$\nbigl(\rho')
\simeq
\psi^{\ast}\bigl(
L\otimes M_1^{-\deg(L)}
\bigr)$.
Take any $b\in B$.
There exists
$\rho\in H^1(Z,U(1))$
such that 
\[
 (L\otimes M_1^{-\deg(L)})_{b}
 \simeq
 \nbigl(\rho)_{b}.
\]
Because
$\nbigl(\rho')_{b}\simeq
\psi^{\ast}(\nbigl(\rho))_{b}$,
we obtain $\rho'=\psi^{\ast}(\rho)$,
and hence
$\psi^{\ast}(L)\simeq
\psi^{\ast}\bigl(
 M_1^{\deg(L)}\otimes \nbigl(\rho)
\bigr)$.
By Corollary \ref{cor;24.10.29.4},
we obtain
$L\simeq M_1^{\deg(L)}\otimes\nbigl(\rho)$.
\hfill\qed

\section{Horizontal deformations of Higgs bundles}
\label{subsection;25.1.22.32}

\subsection{Horizontal deformations of line bundles on smooth spectral curves}

\subsubsection{Basic properties of smooth spectral curves}
Let $X$ be a compact Riemann surface
with the genus $g(X)\geq 2$.
Let $\pi_0:T^{\ast}X\to X$ denote the projection.
Let $\eta\in H^0(T^{\ast}X,\pi_0^{\ast}(K_X))$
denote the tautological section,
i.e., it is the section induced by
the diagonal embedding
$T^{\ast}X\to T^{\ast}X\times_XT^{\ast}X$.
For any positive integer $n$,
we set $\vecA_{X,n}=\bigoplus_{j=1}^nH^0(X,K_X^j)$.
For any $s=(s_j\,|\,j=1,\ldots,n)\in\vecA_{X,n}$,
let $\Sigma_s\subset T^{\ast}X$ denote the $0$-scheme of
the section
$\eta^n+\sum_{j=1}^n(-1)^j\pi_0^{\ast}(s_j)\eta^{n-j}\in
H^0(T^{\ast}X,\pi_0^{\ast}(K_X^n))$.

Let $\vecA'_{X,n}\subset \vecA_{X,n}$ denote the Zariski open subset of
$s\in \vecA_{X,n}$ such that $Z_{s}$ are smooth schemes.
For any $s\in\vecA'_{X,n}$,
let $\pi_{s}:Z_s\to X$
denote the projection.
We recall the following
from \cite{Beauville-Narasimhan-Ramanan}.
\begin{prop}
\label{prop;25.1.22.1}
$Z_s$ is connected.
The genus of $Z_s$
equals
$n^2(g(X)-1)+1$.
\hfill\qed
\end{prop}

It would be instructive to recall the following lemma.
\begin{lem}
\label{lem;24.10.30.1}
There exists a natural isomorphism
$\pi_{s}^{\ast}(K_X^n)\simeq
K_{Z_s}$.
 \end{lem}
\pf
Let $\nbign$ denote 
the normal bundle of $Z_s$ in $T^{\ast}X$.
By the construction,
$\nbign$ equals $\pi_s^{\ast}(K_X^n)$.
By the canonical symplectic structure of $T^{\ast}X$,
$\nbign$ and
the tangent bundle $K_{Z_s}^{-1}$
are mutually dual.
Then, we obtain the desired isomorphism 
$\pi_{s}^{\ast}(K_X^n)\simeq
K_{Z_s}$.
\hfill\qed

\subsubsection{Horizontal line bundles}

We set
\[
 Z_{X,n}=
 \Bigl\{
 (s,\eta)\in \vecA'_{X,n}\times T^{\ast}X\,
 \Big|\,
 \eta^n+\sum_{j=1}^n(-1)^js_{j}\eta^{n-j}=0
 \Bigr\}.
\]

Let $B$ be a multi-disc.
Let $\varphi:B\to \vecA'_{X,n}$ be any holomorphic map.
We set $Z_B=Z_{X,n}\times_{\vecA'_{X,n}}B$.
Let $f_B:Z_B\to B$ denote the projection.
Let $q_{B}:Z_B\to X$ denote the naturally induced morphism.

For any $b\in B$,
we set $Z_b:=Z_{\varphi(b)}$,
and let $\iota_b:Z_b\to Z_B$ denote the inclusion.
For any coherent $\nbigo_{Z_B}$-module $\nbigf$
flat over $f_B^{-1}\nbigo_B$,
we set
$\nbigf_b=\iota_b^{\ast}(\nbigf)$.
We also set
$\deg(\nbigf):=\deg(\nbigf_b)$,
which is independent of $b\in B$.

Let $\Omega^1_{Z_B/B}$
denote the quotient of
$f_B^{\ast}\Omega^1_{B}
\to
\Omega^1_{Z_B}$.
As a straightforward generalization of 
Lemma \ref{lem;24.10.30.1},
we obtain the following natural isomorphism:
\begin{equation}
\Omega^1_{Z_B/B}
\simeq
q_{B}^{\ast}(K_X^n).
\end{equation}

\begin{df}
\label{df;24.10.28.1}
A holomorphic line bundle $L$
on $Z_B$ is called horizontal
if it is horizontal with respect to
$q_{B}^{\ast}(K_X)$.
It is equivalent to the condition that
$L$ is horizontal with respect to $\Omega^1_{Z_B/B}$.
(See Lemma {\rm\ref{lem;24.10.30.2}}.)
\hfill\qed
\end{df}
We reword Proposition \ref{prop;24.10.30.3}
as follows.
\begin{prop}
\label{prop;24.10.30.10}
Let $b\in B$.
For any holomorphic line bundle $N$ on $Z_{b}$,
there exits a horizontal line bundle $N_B$ on $Z_B$
with an isomorphism 
$(N_B)_b\simeq N$.
Such a line bundle is unique up to isomorphisms.
\hfill\qed
\end{prop}

Let $K_X^{1/2}$ be a holomorphic line bundle on $X$
with an isomorphism
$(K_X^{1/2})^2\simeq K_X$.
Note that $\deg\pi_{s}^{\ast}(K_X^{(n-1)/2})=n(n-1)(g(X)-1)$.
The following lemma is clear by the definitions.
\begin{lem}
\label{lem;24.10.30.5}
Let $L$ be a holomorphic line bundle on $Z_B$
with $\deg(L)=n(n-1)(g(X)-1)$.
It is horizontal in the sense of
Definition {\rm\ref{df;24.10.28.1}}
if and only if
there exists
$\rho\in H^1(Z_B,U(1))$ 
such that 
 $L\simeq q_{B}^{\ast}(K_X^{(n-1)/2})\otimes \nbigl(\rho)$,
i.e.,
it is horizontal  
in the sense of  
{\rm\cite{Mochizuki-Asymptotic-Hitchin-metric}.}
\hfill\qed
\end{lem}

\begin{lem}
\label{lem;24.10.29.30}
Let $R$ be any holomorphic line bundle on $X$.
Then,  
$q_{B}^{\ast}(R)$
is horizontal. 
\end{lem}
\pf
Let $K_1$ be a holomorphic line bundle on $X$
with an isomorphism $K_1^{2g(X)-2}\simeq K_X$.
Note that $\deg(K_1)=1$.
There exists
$\rho\in H^1(X,U(1))$
such that
$R\simeq K_1^{\deg(R)}\otimes \nbigl_X(\rho)$,
where $\nbigl_X(\rho)$ denotes the holomorphic line bundle on $X$
underlying the unitary flat line bundle induced by $\rho$.

There exists the natural homomorphism
$q_{B}^{\ast}:H^1(X,U(1))\to H^1(Z_B,U(1))$,
and there exists a natural isomorphism
$q_{B}^{\ast}\nbigl_X(\rho)
\simeq
 \nbigl(q_{B}^{\ast}(\rho))$.
There exists the following natural isomorphism
\[
 q_{B}^{\ast}(R)
 \simeq
 q_{B}^{\ast}(K_1)^{\deg(R)}
 \otimes
 \nbigl(q_{B}^{\ast}(\rho)).
\]
Because $q_{B}^{\ast}(K_1)$ is horizontal,
we obtain that $q_{B}^{\ast}(R)$ is also horizontal.
\hfill\qed

\subsubsection{Pull back by a covering map}

Let $Y$ be a compact Riemann surface
with a holomorphic covering map $\psi:Y\to X$.
It induces $\psi^{\ast}:\vecA_{X,n}\to \vecA_{Y,n}$.
The map $\psi$ induces
a covering map $T^{\ast}Y\to T^{\ast}X$.
For any $s\in\vecA_{X,n}$,
we have
\begin{equation}
\label{eq;24.10.30.4}
Z_{\psi^{\ast}(s)}=Z_s\times_{T^{\ast}X}T^{\ast}Y.
\end{equation}
Then, we obtain the following lemma.
\begin{lem}
We obtain $\psi^{\ast}:\vecA'_{X,n}\to \vecA'_{Y,n}$. 
Moreover,
$Z_{\psi^{\ast}(s)}\to Z_s$ 
is a covering map of Riemann surfaces.
\hfill\qed
\end{lem}

Let $\varphi:B\to \vecA'_{X,n}$ be any holomorphic map.
We set
$Z_{B,\varphi}:=Z_{X,n}\times_{\vecA'_{X,n}}B$.
By using $\psi^{\ast}(\varphi):B\to \vecA'_{Y,n}$,
we also set
$Z_{B,\psi^{\ast}(\varphi)}:=Z_{Y,n}\times_{\vecA'_{Y,n}}B$.
The holomorphic covering map $T^{\ast}Y\to T^{\ast}X$
induces a holomorphic covering map
\[
 \psitilde:
 Z_{B,\psi^{\ast}(\varphi)}\lrarr
 Z_{B,\varphi}.
\]
There exists a natural isomorphism
$\psitilde^{\ast}\Omega^1_{Z_{B,\varphi}/B}
\simeq
 \Omega^1_{Z_{B,\psi^{\ast}(\varphi)}/B}$.
By Lemma \ref{lem;24.10.29.10},
we obtain the following proposition.
\begin{prop}
A holomorphic line bundle $L$ on $Z_{B,\varphi}$ is horizontal
if and only if 
$\psitilde^{\ast}(L)$ on $Z_{B,\psi^{\ast}(\varphi)}$ is horizontal.

\hfill\qed
\end{prop}

\subsection{Horizontal deformations of Higgs bundles with smooth spectral curve}
\label{subsection;25.1.22.31}

\subsubsection{Line bundles and the induced Higgs bundles}

Let $B$ be a multi-disc.
Let $\varphi:B\to \vecA'_{X,n}$ be any holomorphic map.
We set $Z_{B,\varphi}=Z_{X,n}\times_{\vecA'_{X,n}}B$.
Let $f_{B,\varphi}:Z_{B,\varphi}\to B$ denote the projection.
Let $q_{B,\varphi}:Z_{B,\varphi}\to X$ denote the naturally induced morphism.
There exists the induced morphism
$\pi_{B,\varphi}:Z_{B,\varphi}\to B\times X$.
We note that $\pi_{B,\varphi}$ is flat
because $Z_{B,\varphi}$ is defined by
$\eta^n+\sum_{j=1}^n(-1)^j\varphi^{\ast}(s_j)\eta^{n-j}=0$.

Let $L$ be any line bundle on $Z_{B,\varphi}$.
The tautological section $\eta$ induces
$\eta:L\to L\otimes q_{B,\varphi}^{\ast}(K_X)$.
It induces
a morphism
$(\pi_{B,\varphi})_{\ast}(\eta):
 (\pi_{B,\varphi})_{\ast}(L)
 \to
 (\pi_{B,\varphi})_{\ast}(L)
 \otimes
 \Omega^1_{B\times X/B}$,
i.e., 
a relative Higgs field of $(\pi_{B,\varphi})_{\ast}(L)$.
In this way,
we obtain a relative Higgs bundle
$(\pi_{B,\varphi})_{\ast}(L,\eta):=
\bigl((\pi_{B,\varphi})_{\ast}(L),
 (\pi_{B,\varphi})_{\ast}(\eta)
\bigr)$.
The restrictions of 
$(\pi_{B,\varphi})_{\ast}(L)_{|\{b\}\times X}$
and
$(\pi_{B,\varphi})_{\ast}(\eta)_{|\{b\}\times X}$
to $\{b\}\times X$
are denoted by
$(\pi_{B,\varphi})_{\ast}(L)_b$
and
$(\pi_{B,\varphi})_{\ast}(\eta)_b$.

\begin{lem}
We have 
$\deg(L_b)-\deg\bigl(
 (\pi_{B,\varphi})_{\ast}(L)_{b}
 \bigr)
=n(n-1)(g(X)-1)$
for any $b\in B$.
\end{lem}
\pf
By the Riemann-Roch theorem
and Proposition \ref{prop;25.1.22.1},
we obtain
\[
 \deg(L_b)+1-\bigl(
 n^2(g(X)-1)+1
 \bigr)
 =\deg\bigl(
 (\pi_{B,\varphi})_{\ast}(L)_b
 \bigr)
 +n(1-g(X)).
\]
Then, we obtain the claim of the lemma.
\hfill\qed

\subsubsection{Horizontal deformations}
\label{subsection;25.1.22.30}

Let $(\nbige,\theta)$ be a relative Higgs bundle of rank $n$
on $B\times X/B$
such that $\Sigma_{\theta}=Z_{B,\varphi}$.
According to \cite{Beauville-Narasimhan-Ramanan},
there exists a line bundle $L_{\nbige}$ on $Z_{B,\varphi}$
such that
$(\nbige,\theta)\simeq (\pi_{B,\varphi})_{\ast}(L_{\nbige},\eta)$.

\begin{df}
\label{df;24.10.29.11}
$(\nbige,\theta)$ is called horizontal
if $L_{\nbige}$ is horizontal
in the sense of  Definition {\rm\ref{df;24.10.28.1}}.
\hfill\qed
\end{df}

The restriction of $(\nbige,\theta)$
to $\{b\}\times X$
is denoted by $(\nbige_b,\theta_b)$.
We obtain the following lemma from Lemma \ref{lem;24.10.30.5}.
\begin{lem}
If $\deg(\nbige_b)=0$,
then $\nbige$ is horizontal in the sense of
Definition {\rm\ref{df;24.10.29.11}}
if and only if 
it is horizontal in the sense of
{\rm\cite{Mochizuki-Asymptotic-Hitchin-metric}}.
\hfill\qed
\end{lem}

We obtain the following proposition from 
Proposition \ref{prop;24.10.30.10}.
\begin{prop}
\label{prop;24.10.29.20}
Let $b\in B$.
For any Higgs bundle $(E,\theta_E)$ on $X$
with $\Sigma_{\theta_E}=Z_{b}$,
there exists a relative Higgs bundle
$(\nbige,\theta)$ on $B\times X$
such that 
(i) $\Sigma_{\theta}=Z_{B,\varphi}$,
(ii) $(\nbige,\theta)$ is horizontal,
(iii) $(\nbige_{b},\theta_{b})\simeq (E,\theta_E)$.
Such a relative Higgs bundle is unique up to isomorphisms. 
\hfill\qed
\end{prop}

We obtain the following lemma from
Lemma \ref{lem;24.10.29.31}
and Lemma \ref{lem;24.10.29.30}.
\begin{lem}
\label{lem;24.10.29.32}
Let $(\nbige,\theta)$ be a horizontal relative Higgs bundle
on $B\times X$ over $B$.
Let $L$ be a line bundle on $X$.
The pull back of $L$  to $B\times X$ is also denoted by $L$.
Then, $(\nbige\otimes L,\theta)$ is also horizontal.
\hfill\qed
\end{lem}

Let $\psi:Y\to X$ be a covering map of compact Riemann surfaces.
We obtain $\psi^{\ast}(\varphi):B\to \vecA'_{Y,n}$
and $Z_{B,\psi^{\ast}(\varphi)}=Z_{Y,n}\times_{\vecA'_{Y,n}}B$.
There exists the covering map
$\psitilde:Z_{B,\psi^{\ast}(\varphi)}\to Z_{B,\varphi}$.
Let $\eta_Y$ denote the tautological $1$-form on $T^{\ast}Y$.
For a line bundle $L_Y$ on $Z_{B,\psi^{\ast}(\varphi)}$,
we obtain the relative Higgs bundle
$(\pi_{B,\psi^{\ast}(\varphi)})_{\ast}(L_Y,\eta_Y)$
on $B\times Y$ over $B$.

\begin{prop}
A relative Higgs bundle $(\nbige,\theta)$ on $B\times X$ over $B$
with $\Sigma_{\theta}=Z_{B,\varphi}$
is horizontal
if and only if
the relative Higgs bundle
$(\id_B\times\psi)^{\ast}(\nbige,\theta)$
on $B\times Y$ over $B$
is horizontal.
\end{prop}
\pf
Let $L_{\nbige}$ be a holomorphic line bundle on $Z_{B,\varphi}$
with an isomorphism
$(\pi_{B,\varphi})_{\ast}(L_{\nbige},\eta)\simeq (\nbige,\theta)$.
Because $\id\times \psi$ is a holomorphic covering map,
there exists the following isomorphism
by the standard base change theorem:
\[
 (\pi_{B,\psi^{\ast}(\varphi)})_{\ast}\bigl(
 \psitilde^{\ast}(L_{\nbige}),\eta_Y
 \bigr)
 \simeq
 (\id_B\times \psi)^{\ast}
 (\pi_{B,\varphi})_{\ast}\bigl(
 L_{\nbige},\eta
 \bigr)
 \simeq
 (\id_B\times \psi)^{\ast}
 (\nbige,\theta).
\]
We obtain 
$\Sigma_{\psi^{\ast}(\theta)}=Z_{B,\psi^{\ast}(\varphi)}$.
By Lemma \ref{lem;24.10.29.10},
$L_{\nbige}$ is horizontal if and only if $\psitilde^{\ast}(L_{\nbige})$.
Hence, we obtain the claim of the proposition.
\hfill\qed

\section{The comparison of Hitchin metric and the semi-flat metric}

\subsection{The moduli space of stable Higgs bundles}

For $(n,d)\in\seisuu_{>0}\times\seisuu$,
let $\nbigm_H(X,n,d)$ denote the moduli space of
stable Higgs bundles of rank $n$ and degree $d$ on $X$.

\subsubsection{Deformation theory}

Let $(E,\theta)\in\nbigm_H(X,n,d)$.
Let $\Def(E,\theta)$ denote the complex of sheaves
$\End(E)
 \stackrel{\ad\theta}{\lrarr}
 \End(E)\otimes K_X$,
where $\End(E)$ sits at the degree $0$.
There exists the natural isomorphism:
\[
 T_{(E,\theta)}\nbigm_H(X,n,d)
 \simeq
 H^1(X,\Def(E,\theta)).
\]

Let $A^j(\End(E))$ denote the space pf
$\End(E)$-valued $j$-forms.
Together with the differential
\[
 \delbar_E+\ad\theta:
A^j(\End(E))
\to A^{j+1}(\End(E)),
\]
we obtain the complex
$A^{\bullet}(\End(E),\ad\theta)$.
It is also obtained as the space of global sections
of the Dolbeault resolution of $\Def(E,\theta)$.
There exists the natural isomorphism
\[
 H^1(\Def(E,\theta))
 \simeq
 H^1\bigl(A^{\bullet}(\End(E,\theta))\bigr).
\]

\subsubsection{Hitchin metric}
\label{subsection;25.1.22.10}

Let $(E,\theta)\in\nbigm_H(X,n,d)$.
For any Hermitian metric $h$ of $E$,
let $R(h)$ denote the curvature of
the Chern connection of $(E,h)$.
Let $\theta^{\dagger}_h$ denote the adjoint of $\theta$
with respect to $h$.
We set
$F(h)=R(h)+[\theta,\theta^{\dagger}_h]$.
A Hermitian metric $h$ induces a Hermitian metric of $\End(E)$
denoted by $h_{\End(E)}$.
We have
\[
 F(h_{\End(E)})
=\ad F(h_{\End(E)}).
\]

We say $h$ is a Hermitian-Einstein metric of $(E,\theta)$
if the trace-free part of $F(h)$ is $0$,
i.e.,
\[
 F(h)-\frac{1}{n}\tr F(h)\id_E=0.
\]
The induced Hermitian metric $h_{\End(E)}$
is a harmonic metric of
the Higgs bundle $(\End(E),\ad\theta)$,
i.e.,
$F(h_{\End(E)})=0$.

\begin{thm}[Hitchin-Simpson,
\cite{Hitchin-self-duality}, \cite{s1}]
 \mbox{{}}\label{thm;24.10.28.1}
\begin{itemize}
 \item $(E,\theta)$ has a Hermitian-Einstein metric $h$.
 \item Let
       $h_j$ $(j=1,2)$ be Hermitian-Einstein metrics of $(E,\theta)$.
       Let $\varphi:X\to\real_{>0}$ be the $C^{\infty}$-function
       determined by
       $\det(h_1)=\varphi^n\det(h_2)$.
       Then,
       we obtain $h_1=\varphi h_2$.
       In particular,
       we obtain
       $(h_1)_{\End(E)}=(h_2)_{\End(E)}$.
       \hfill\qed
\end{itemize} 
 \end{thm}
In the following,
let $h_{\HE}$ denote a Hermitian-Einstein metric of $(E,\theta)$,
and let $h_{\HE,\End(E)}$ denote the induced harmonic metric
of $(\End(E),\ad\theta)$
in Theorem \ref{thm;24.10.28.1}.

For any $\End(E)$-valued $1$-forms $\tau_i$ $(i=1,2)$,
we set
\[
 (\tau_1,\tau_2)_{h_{\HE,\End(E)},L^2}
 =2\sqrt{-1}
 \int_X
 \tr\Bigl(
 \tau_1^{1,0}
 \wedge
 (\tau_2^{1,0})^{\dagger}_{h_{\HE}}
-\tau_1^{0,1}
 \wedge
 (\tau_2^{0,1})^{\dagger}_{h_{\HE}}
\Bigr).
\]
It is independent of the choice of $h_{\HE}$.

Let $\nbigh^1(\End(E),\theta)$
denote the space of harmonic $1$-forms of
$(\End(E),\ad(\theta),h_{\End(E)})$.
There exists the standard isomorphism
\[
 \nbigh^1\bigl(\End(E),\ad\theta\bigr)
 \simeq
 H^1(X,\Def(E,\theta)),
 \quad
 \tau\longmapsto[\tau].
\]
The Hitchin metric $g^{X,n,d}_{H}$ at $(E,\theta)$
is defined as
\[
 g^{X,n,d}_{H|(E,\theta)}([\tau_1],[\tau_2])
=\bigl(\tau_1,\tau_2\bigr)_{L^2,h_{\HE,\End(E)}}.
\]

\subsubsection{Complex symplectic structure}

For $\tau_1,\tau_2\in A^1(\End(E),\ad\theta)$,
we define
\[
 \omega^{X,n,d}_{H|(E,\theta)}(\tau_1,\tau_2)
=-\int_X\Tr(\tau_1\wedge\tau_2).
\]
This induces an alternative form on
$H^1\bigl(A^{\bullet}(\End(E),\ad\theta)\bigr)$,
which we denote by the same notation.
According to Hitchin \cite{Hitchin-self-duality},
the alternative forms
$\omega^{X,n,d}_{H|(E,\theta)}$ $((E,\theta)\in\nbigm_H(X,n,d))$
define a complex symplectic structure
$\omega^{X,n,d}_{H}$ of
$\nbigm_H(X,n,d)$.

\subsubsection{Twist by a line bundle}

Let $L$ be a holomorphic line bundle of degree $e$ on $X$.
By taking the tensor product with $L$,
we obtain the isomorphism
$F_L:\nbigm_H(X,n,d)\simeq \nbigm_H(X,n,d+ne)$,
i.e.,
$F_L(E,\theta)=(E\otimes L,\theta)$.
The following lemma is clear by the constructions.
\begin{lem}
\label{lem;24.10.29.40}
 We have
 $F_L^{\ast}\bigl(
 g^{X,n,d+ne}_{H}
 \bigr)
 =g^{X,n,d}_{H}$
and 
 $F_L^{\ast}\bigl(
 \omega^{X,n,d+ne}_{H}
 \bigr)
 =\omega^{X,n,d}_{H}$.
\hfill\qed
\end{lem}

\subsection{The open subsets of Higgs bundles with smooth spectral curve}

\subsubsection{Hitchin fibration}
\label{subsection;25.1.22.12}

For any $(E,\theta)\in\nbigm_{H}(X,n,d)$,
we obtain the section
\[
 \det(\eta\id_E-\theta)
 =\eta^n
 +\sum_{j=1}^n(-1)^ja_j(E,\theta)\eta^{n-j}
 \in
 H^0(T^{\ast}X,\pi_0^{\ast}(K_X^n)).
\]
We obtain the morphism
$\Phi_{n,d}:\nbigm_H(X,n,d)\to\vecA_{X,n}$
by
$\Phi_{n,d}(E,\theta)=(a_j(E,\theta)\,|\,j=1,\ldots,n)$,
which is called the Hitchin fibration.
We set
$\nbigm_H'(X,n,d)=\Phi_{n,d}^{-1}(\vecA'_{X,n})$.
We also set
$\Sigma_{\theta}=\Sigma_{\Phi_{n,d}(E,\theta)}$.
We recall the following proposition
from \cite{Beauville-Narasimhan-Ramanan} and
\cite{Hitchin-self-duality}.
\begin{prop}
\mbox{{}}
\begin{itemize}
 \item The induced morphism
$\nbigm_H'(X,n,d)\to\vecA'_{X,n}$ is smooth.
 \item
For any $s\in\vecA'_{X,n}$,
      $\Phi_{n,d}^{-1}(s)$ is isomorphic to
$\Pic_{N(n,d)}(\Sigma_{s})$, where
$N(n,d)=d+n(n-1)(g(X)-1)$.
\hfill\qed
\end{itemize}
\end{prop}

We define $(T_{(E,\theta)}\nbigm_{H}(X,n,d))^{\ver}$
as the image of
\[
 T_{(E,\theta)}\Phi_{n,d}^{-1}\bigl(
 \Phi_{n,d}(E,\theta)
 \bigr)
 \lrarr
 T_{(E,\theta)}\nbigm'_{H}(X,n,d).
\]
It is naturally isomorphic to
$H^1(\Sigma_{\theta},\nbigo_{\Sigma_{\theta}})$.
The following lemma is also due to Hitchin \cite{Hitchin-self-duality}.
\begin{lem}
$(T_{(E,\theta)}\nbigm_{H}(X,n,d))^{\ver}$
is Lagrangian with respect to
$\omega^{X,n,d}_{H}$.
\hfill\qed
\end{lem}

\subsubsection{Covering map}

Suppose that $g(X)\geq 2$.
Let $\psi:Y\to X$ be a covering map of compact Riemann surfaces.
We assume that it is a Galois covering,
i.e.,
there exists a finite group $G$ acting freely on $Y$
such that $Y/G\simeq X$.
Let $e$ denote the degree of $\psi$.
\begin{rem}
For any $e\in\seisuu_{>0}$,
there exists a Galois covering
such that $G$ is a cyclic group $\seisuu/e\seisuu$. 
\hfill\qed
\end{rem}

Let $(E,\theta)$ be a stable Higgs bundle on $X$.
We obtain the induced Higgs bundle
$\psi^{\ast}(E,\theta)$ on $Y$,
which is naturally $G$-equivariant.
Because
$\Sigma_{\psi^{\ast}\theta}=\Sigma_{\theta}\times_XY$,
we obtain the following lemma.

\begin{lem}
\label{lem;24.10.30.12}
This induces a morphism
$F_{\psi}:\nbigm'_{H}(X,n,d)\to \nbigm'_{H}(Y,n,ed)$.
\hfill\qed
\end{lem}

Let $(E,\theta)\in\nbigm'_{H}(X,n,d)$.
Because
$\Def(\psi^{\ast}E,\psi^{\ast}\theta)
=\psi^{\ast}(\Def(E,\theta))$
is naturally $G$-equivariant,
the space
$H^1(Y,\Def(\psi^{\ast}E,\psi^{\ast}\theta))$
is naturally a $G$-representation.

\begin{lem}
The natural morphism
$H^1(X,\Def(E,\theta))\to
H^1(Y,\Def(\psi^{\ast}E,\psi^{\ast}\theta))$
is injective,
 and the image is
the $G$-invariant part of 
$H^1(Y,\Def(\psi^{\ast}E,\psi^{\ast}\theta))$.
\end{lem}
\pf
Note that $\psi_{\ast}(\nbigo_{Y})$
is equipped with a naturally defined $G$-action,
and that there exists the natural morphism
$\nbigo_X\to \psi_{\ast}(\nbigo_{Y})$.
Moreover,
there exists the decomposition
of locally free $\nbigo_X$-modules
\[
 \psi_{\ast}(\nbigo_{Y})
 =\nbigo_X
 \oplus
 \nbigf
\]
equipped with $G$-actions.
Each fiber $\nbigf_{|P}$ $(P\in X)$
does not contain a trivial $G$-representation.
We obtain
\[
\psi_{\ast}\psi^{\ast}\Def(E,\theta)
=\Def(E,\theta)\otimes
\psi_{\ast}(\nbigo_Y)
=\Def(E,\theta)\oplus
\Def(E,\theta)\otimes\nbigf.
\]
Then, the claim is clear.
\hfill\qed

\vspace{.1in}
We may also obtain the lemma by considering harmonic $1$-forms.
Let $h_{\HE}$ be a Hermitian-Einstein metric of $(E,\theta)$.
Then, $\psi^{\ast}(h_{\HE})$ is also a Hermitian-Einstein metric of
$\psi^{\ast}(E,\theta)$.
The following lemma is clear.
\begin{lem}
The pull back induces
$\psi^{\ast}:\nbigh^1(\End(E),\ad\theta)
\to \nbigh^1(\End(\psi^{\ast}(E)),\ad\psi^{\ast}\theta)$.
It is injective,
and the image equals the $G$-invariant part of
$\nbigh^1(\End(\psi^{\ast}(E)),\ad\psi^{\ast}\theta)$.
\hfill\qed
\end{lem}

\begin{prop}
\label{prop;24.10.29.41}
We have
\[
 (T_{(E,\theta)}F_{\psi})^{\ast}(g^{Y,n,de}_H)_{|(E,\theta)}
 =e\cdot g^{X,n,d}_{H|(E,\theta)},
\quad\quad
 (T_{(E,\theta)}F_{\psi})^{\ast}
 (\omega^{Y,n,de}_H)_{|(E,\theta)}
=e\cdot \omega^{X,n,d}_{H|(E,\theta)}.
\] 
\end{prop}
\pf
Let $\tau_i$ $(i=1,2)$ be harmonic $1$-forms of
$(\End(E),\ad\theta,h_{\HE,\End(E)})$.
We obtain
\begin{multline}
 2\sqrt{-1}\int_{Y}
 \Tr\Bigl(
 \psi^{\ast}(\tau_1)^{1,0}
 \bigl(
 \psi^{\ast}(\tau_2)^{1,0}
 \bigr)^{\dagger}_{\psi^{\ast}(h_{\HE})}
-
 \psi^{\ast}(\tau_1)^{0,1}
 \bigl(
 \psi^{\ast}(\tau_2)^{0,1}
 \bigr)^{\dagger}_{\psi^{\ast}(h_{\HE})}
 \Bigr)
 = \\
 e\times 2\sqrt{-1}\int_X
  \Tr\Bigl(
  \tau_1^{1,0}
  \bigl(
  \tau_2^{1,0}
  \bigr)^{\dagger}_{h_{\HE}}
-
 \tau_1^{0,1}
 \bigl(
 \tau_2^{0,1}
 \bigr)^{\dagger}_{h_{\HE}}
 \Bigr).
\end{multline}
Hence, we obtain 
$(T_{(E,\theta)}F_{\psi})^{\ast}(g^{Y,n,de}_H)_{(E,\theta)}
 =e\cdot g^{X,n,d}_{H|(E,\theta)}$.
The other is obtained similarly.
\hfill\qed

\subsection{Semi-flat metrics in the non-zero degree case}

\subsubsection{Locally principal bundles and integrable connections}
\label{subsection;24.10.30.11}

For any $s\in \vecA'_{X,n}$,
let $U_{s}$ be a neighbourhood of $s$ in $\vecA'_{X,n}$,
which is isomorphic to a multi-disc.
For any $s'\in U_{s}$,
we set
$T(s'):=H^1(\Sigma_{s'},U(1))$.
We can regard $T(s')$ as
the moduli space of flat line bundles on $\Sigma_{s'}$,
which naturally acts on 
$\Pic_{N(n,d)}(\Sigma_{s'})$
by the tensor product of the underlying holomorphic line bundles.
The action is transitive and free.
A path connecting $s$ and $s'$
induces an isomorphism $T(s)\simeq T(s')$.
Hence,
we obtain the right $T(s)$-action on $\Phi_{n,d}^{-1}(U_s)$:
\[
 \rho:\Phi_{n,d}^{-1}(U_s)\times T(s)\lrarr
 \Phi_{n,d}^{-1}(U_s),
\]
with which
$\Phi_{n,d}^{-1}(U_{s})\to U_{s}$
is a principal $T(s)$-bundle.

Let $(E,\theta_E)\in\Phi_{n,d}^{-1}(s)$.
According to Proposition \ref{prop;24.10.29.20},
there exists a relative Higgs bundle
$(\nbige,\theta)$ on $U_s\times X$ over $U_s$
such that
(i) $\Sigma_{\theta}=U_s\times_{\vecA'_{X,n}}Z_{X,n}$,
(ii) $(\nbige,\theta)$ is horizontal,
(iii) $(\nbige,\theta)_{|s}\simeq (E,\theta_E)$.
It induces a section
$H_{(E,\theta)}:U_s\to \nbigm'_H(X,n,d)\times_{\vecA'_{X,n}}U_s$
of the Hitchin fibration.
By the construction,
we have
\[
 H_{\rho((E,\theta),a)}(s')=
 \rho\bigl(
 H_{(E,\theta)}(s'),a
 \bigr)
 \quad
 (s'\in U_s,\,\,a\in T(s)).
\]
Therefore,
we obtain the unique integrable connection
of the $T(s)$-principal bundle
$\Phi_{n,d}^{-1}(U_s)$ over $U_s$
such that
$H_{(E,\theta)}$
$((E,\theta)\in\Phi_{n,d}^{-1}(s))$
are horizontal sections.

\subsubsection{Lagrangian property of the horizontal sections}

We define
\[
 \bigl(
 T_{(E,\theta)}\nbigm'_H(X,n,d)
 \bigr)^{\hor}
 \subset
 T_{(E,\theta)}\nbigm'_H(X,n,d)
\]
as the image of
$T_sH_{(E,\theta)}:
T_s\vecA'_{X,n}\to T_{(E,\theta)}\nbigm'_H(X,n,d)$.
It is the horizontal subspace
with respect to the integrable connection
in \S\ref{subsection;24.10.30.11}.

\begin{prop}
Each $(T_{(E,\theta)}\nbigm'_H(X,n,d))^{\hor}$
is Lagrangian with respect to
$\omega^{X,n,d}_H$.
\end{prop}
\pf
Let $\psi:Y\to X$ be
a Galois covering map of a compact Riemann surfaces
such that the degree $e$ of $\psi$ is a multiple of $n$.
Let $L$ be a holomorphic line bundle on $Y$
with $\deg(L)=(ed)/n$.

For any $(E,\theta)\in\nbigm_{H}'(X,n,d)$,
we obtain
$F_{\psi,L}(E,\theta)=
\psi^{\ast}(E,\theta)\otimes L^{-1}
 \in
 \nbigm_{H}(Y,n,0)$.
This induces a morphism
\[
F_{\psi,L}:\nbigm'_{H}(X,n,d)\to \nbigm'_{H}(Y,n,0).
\]
By Lemma \ref{lem;24.10.29.32},
the derivative
$T_{(E,\theta)}F_{\psi,L}:
T_{(E,\theta)}\nbigm'_{H}(X,n,d)\to
T_{F_{\psi,L}(E,\theta)}\nbigm'_{H}(Y,n,0)$
of $F_{\psi,L}$ at $(E,\theta)$
induces
\[
\bigl(
T_{(E,\theta)}\nbigm'_{H}(X,n,d)
\bigr)^{\hor}
\lrarr
\bigl(
T_{F_{\psi,L}(E,\theta)}\nbigm'_{H}(Y,n,0)
\bigr)^{\hor}.
\]
According to \cite[Corollary 3.33]{Mochizuki-Asymptotic-Hitchin-metric},
$\bigl(
T_{F_{\psi,L}(E,\theta)}\nbigm'_{H}(Y,n,0)
\bigr)^{\hor}$
is Lagrangian with respect to
$\omega^{Y,n,0}_H$.
By Lemma \ref{lem;24.10.29.40}
and Proposition \ref{prop;24.10.29.41},
the restriction of
$\omega^{X,n,d}_H$
to
$(T_{(E,\theta)}\nbigm_{H}(X,n,d))^{\hor}$ is $0$.
Because
\[
2\dim
(T_{(E,\theta)}\nbigm_{H}(X,n,d))^{\hor}
=\dim
T_{(E,\theta)}\nbigm_{H}(X,n,d),
\]
$(T_{(E,\theta)}\nbigm_{H}(X,n,d))^{\hor}$
is Lagrangian.
\hfill\qed

\subsubsection{The tangent spaces
and the cohomology group of the spectral curves}

Let $(E,\theta)\in\nbigm'_H(X,n,d)$.
We obtain the decomposition
into the vertical direction and the horizontal direction:
\[
T_{(E,\theta)}\nbigm'_{H}(X,n,d) 
=(T_{(E,\theta)}\nbigm'_{H}(X,n,d))^{\ver}
\oplus
(T_{(E,\theta)}\nbigm'_{H}(X,n,d))^{\hor}.
\]

Let $\pi_{\theta}:\Sigma_{\theta}\to X$
denote the projection.
Let $\nbigh^j(\Def(E,\theta))$
denote the $j$-th cohomology sheaves.
As explained in \cite{Mochizuki-Asymptotic-Hitchin-metric},
there exist the following natural isomorphisms
\[
 (\pi_{\theta})_{\ast}
 (\nbigo_{\Sigma_{\theta}})
 \simeq
 \nbigh^0\bigl(
 \Def(E,\theta)
 \bigr),
\quad\quad
(\pi_{\theta})_{\ast}(K_{\Sigma_{\theta}})
\simeq
\nbigh^1\bigl(
 \Def(E,\theta)
 \bigr).
\]
We obtain the following exact sequence:
\[
\begin{CD}
 0 @>>>
 H^1(\Sigma_{\theta},\nbigo_{\Sigma_{\theta}})
 @>{a_1}>>
 H^1(X,\Def(E,\theta))
 @>{a_2}>>
 H^0(\Sigma_{\theta},K_{\Sigma_{\theta}})
 @>>>
 0.
\end{CD}
\]
The image of $a_1$
is identified with
$(T_{(E,\theta)}\nbigm'_H(X,n,d))^{\ver}
\subset
T_{(E,\theta)}\nbigm'_H(X,n,d)$,
i.e., we obtain the isomorphism
\[
 \iota^{\ver}_{(E,\theta)}:
 H^1(\Sigma_{\theta},\nbigo_{\Sigma_{\theta}})
 \simeq
 (T_{(E,\theta)}\nbigm'_H(X,n,d))^{\ver}.
\]
We obtain
the isomorphism 
$\iota^{\hor}_{(E,\theta)}:
H^0(\Sigma_{\theta},K_{\Sigma_{\theta}})
\simeq
(T_{(E,\theta)}\nbigm'_H(X,n,d))^{\hor}$
as the composition of the following morphisms:
\[
H^0(\Sigma_{\theta},K_{\Sigma_{\theta}})
\simeq
\left(
T_{(E,\theta)}\nbigm'_H(X,n,d)\Big/
(T_{(E,\theta)}\nbigm'_H(X,n,d))^{\ver}
\right)
\simeq
(T_{(E,\theta)}\nbigm'_H(X,n,d))^{\hor}.
\]

Let $\Harm^1(\Sigma_{\theta})$
denote the space of harmonic $1$-forms
on $\Sigma_{\theta}$.
We obtain
\[
T_{(E,\theta)}\nbigm'_H(X,n,d)
\simeq
H^1(\Sigma_{\theta},\nbigo_{\Sigma_{\theta}})
\oplus
H^0(\Sigma_{\theta},K_{\Sigma_{\theta}})
\simeq
\Harm^1(\Sigma_{\theta}).
\]

\subsubsection{Semi-flat metric}

We recall that
$\Harm^1(\Sigma_{\theta})$
is equipped with the Hermitian metric
defined as follows:
\[
 (\sigma_1,\sigma_2)_{L^2,\Sigma_{\theta}}
 =2\sqrt{-1}
 \int_{\Sigma_{\theta}}
 \Bigl(
 \sigma_1^{1,0}\wedge\overline{\sigma_2^{1,0}}
-\sigma_1^{0,1}\wedge\overline{\sigma_2^{0,1}}
 \Bigr).
\]
Let $g^{X,n,d}_{\semiflat|(E,\theta)}$
be the Hermitian metric
of $T_{(E,\theta)}\nbigm'_{H}(X,n,d)$
induced by
$\bigl(\cdot,\cdot\bigr)_{L^2,\Sigma_{\theta}}$
on $\Harm^1(\Sigma_{\theta})$.

\begin{prop}
\label{prop;25.1.22.100}
$g^{X,n,d}_{\semiflat}=
\bigl\{
g^{X,n,d}_{\semiflat|(E,\theta)}\,\big|\,
(E,\theta)\in\nbigm'_{H}(X,n,d)
\bigr\}$
is a hyperk\"ahler metric.
\end{prop}
\pf
It is enough to check the claim locally.
Let $s\in \vecA'_{X,n}$.
Let $U_s$ be a small neighbourhood.
Let $f:Z_{U_s}\to U_s$
be the induced morphism.
Let $\Theta_{Z_{U_s}/U_s}$ denote the relative tangent sheaf
of $Z_{U_s}$ over $U_s$.
We obtain the locally free $\nbigo_{U_s}$-module
$f_{\ast}(\Theta_{Z_{U_s}/U_s})$.
Let $V$ be the corresponding holomorphic vector bundle
on $U_s$.
By taking the exponential,
the fibers
$V_{|s'}$ $(s'\in U_s)$ acts on $f^{-1}(s')$.
We obtain the lattices $\Lambda_{s'}\subset V_{|s'}$
as the kernel of
$V_{s'}\to \Aut(f^{-1}(s'))$.
By taking the quotient,
we obtain
the family of complex torus
$\nbiga_{U_s}\to U_s$.
According to Freed \cite{Freed},
$\nbiga_{U_s}$ is naturally equipped with
the hyperkahler metric.

Let $H:U_s\to \Phi_H^{-1}(U_s)$
be any horizontal section.
By using the $T(s)$-action,
we obtain the diffeomorphism
$\nbiga_{U_s}\simeq \Phi_H^{-1}(U_s)$.
By the construction,
it is an isometry.
Hence, $g^{X,n,d}_{\semiflat}$ is a hyperk\"ahler metric.
\hfill\qed

\subsubsection{Twist by a line bundle}

Let $L$ be a holomorphic line bundle on $X$
of degree $e$.
Let $F_L:\nbigm'_H(X,n,d)\simeq \nbigm'_H(X,n,d+ne)$
defined by
$F_L(E,\theta)=(E\otimes L,\theta)$.
The following lemma is clear by the construction.
\begin{lem}
$F_L^{\ast}(g^{X,n,d+ne}_{\semiflat})
=g^{X,n,d}_{\semiflat}$.
\hfill\qed
\end{lem}

\subsubsection{Covering map}

Let $\psi:Y\to X$ be a Galois covering map of compact Riemann surfaces.
Let $e$ denote the degree of $\psi$.
We obtain
$F_{\psi}:\nbigm'_H(X,n,d)\to \nbigm'_H(Y,n,de)$
as in Lemma \ref{lem;24.10.30.12}.

\begin{prop}
\label{prop;24.10.30.30}
We have
\begin{equation}
\label{eq;24.10.30.20}
 (T_{(E,\theta)}F_{\psi})^{\ast}
 (g^{Y,n,de}_{\semiflat})_{|(E,\theta)}
 =e\cdot
 g^{X,n,d}_{\semiflat|(E,\theta)}.
\end{equation}
\end{prop}
\pf
The linear map 
$(T_{(E,\theta)}F_{\psi})^{\ast}$
preserves the decomposition into
the vertical part and the horizontal part.
Let $\psitilde:\Sigma_{\psi^{\ast}\theta}\to\Sigma_{\theta}$
denote the induced covering map.
We obtain the following commutative diagram:
\[
\begin{CD}
\Harm^1(\Sigma_{\theta})=
H^1(\Sigma_{\theta},\nbigo_{\Sigma_{\theta}})
\oplus
H^0(\Sigma_{\theta},K_{\Sigma_{\theta}})
@>{\simeq}>{\iota^{\ver}_{(E,\theta)}\oplus \iota^{\hor}_{(E,\theta)}}> 
T_{(E,\theta)}\nbigm_H'(X,n,d)\\
@V{\psitilde^{\ast}}VV @V{T_{(E,\theta)}F_{\psi}}VV \\ 
\Harm^1(\Sigma_{\psi^{\ast}\theta})
=H^1(\Sigma_{\psi^{\ast}\theta},\nbigo_{\Sigma_{\psi^{\ast}\theta}})
\oplus
H^0(\Sigma_{\psi^{\ast}\theta},K_{\Sigma_{\psi^{\ast}\theta}})
 @>{\simeq}>{\iota^{\ver}_{\psi^{\ast}(E,\theta)}
  \oplus \iota^{\hor}_{\psi^{\ast}(E,\theta)}}> 
T_{\psi^{\ast}(E,\theta)}\nbigm_H'(Y,n,de).
\end{CD}
\]
Let $\sigma_1,\sigma_2\in\Harm(\Sigma_{\theta})$.
Because the degree of $\psitilde$ is $e$,
we obtain
\[
 2\sqrt{-1}
 \int_{\Sigma_{\psi^{\ast}\theta}}
 \Bigl(
 \psitilde^{\ast}\sigma_1^{1,0}
 \wedge\overline{\psitilde^{\ast}\sigma^{1,0}_2}
- \psitilde^{\ast}\sigma_1^{0,1}
 \wedge\overline{\psitilde^{\ast}\sigma^{0,1}_2}
 \Bigr)
=e\times(2\sqrt{-1})
 \int_{\Sigma_{\theta}}
 \Bigl(
 \sigma_1^{1,0}
 \wedge\overline{\sigma_2^{1,0}}
-\sigma_1^{0,1}
 \wedge\overline{\sigma_2^{0,1}}
 \Bigr).
\]
Hence, we obtain (\ref{eq;24.10.30.20}).
\hfill\qed

\subsection{Comparison of the Hitchin metric and the semi-flat metric}

Let $s$ be the automorphism of
$T\nbigm'_{H}(X,n,d)$
determined by
$g^{X,n,d}_H=g^{X,n,d}_{\semiflat}\cdot s^{X,n,d}$.

\begin{thm}
\label{thm;25.1.22.40}
For any $(E,\theta)\in \nbigm'_{H}(X,n,d)$,
there exists $\epsilon>0$
such that
\[
\bigl|
 (s^{X,n,d}-\id)_{|(E,t\theta)}
 \bigr|_{g^{X,n,d}_{\semiflat}}
=O\bigl(e^{-\epsilon t}\bigr).
\]
\end{thm}
\pf
Let $\psi:Y\to X$ be a Galois covering
of compact Riemann surfaces
such that the degree $e$ of $\psi$ is a multiple of $n$.
Let $L$ be a line bundle on $Y$
of degree $nd/e$.
We obtain the morphism
\[
 F_{\psi,L}:
 \nbigm'_{H}(X,n,d)
 \lrarr
 \nbigm'_{H}(Y,n,0)
\]
by
$F_{\psi,L}(E,\theta)
=\psi^{\ast}(E,\theta)\otimes L^{-1}$.
We obtain the derivatives
\[
 T_{(E,t\theta)}F_{\psi,L}:
 T_{(E,t\theta)}
 \nbigm'_{H}(X,n,d)
 \lrarr
 T_{F_{\psi,L}(E,t\theta)}
 \nbigm'_{H}(Y,n,0).
\]
By Proposition \ref{prop;24.10.29.41} and
Proposition \ref{prop;24.10.30.30},
we obtain
\[
 (T_{(E,t\theta)}F_{\psi,L})^{\ast}
 (g^{Y,n,0}_{H})_{|(E,t\theta)}
 =e\cdot g^{X,n,d}_{H|(E,t\theta)},
\quad\quad
  (T_{(E,t\theta)}F_{\psi,L})^{\ast}
 (g^{Y,n,0}_{\semiflat})_{|(E,t\theta)}
 =e\cdot
 g^{X,n,d}_{\semiflat|(E,t\theta)}.
\]
By \cite{Mochizuki-Asymptotic-Hitchin-metric},
there exists $\epsilon>0$ such that 
\[
\Bigl|
 g^{Y,n,0}_{\semiflat|F_{L,\psi}(E,t\theta)}
-g^{Y,n,0}_{H|F_{L,\psi}(E,t\theta)}
\Bigr|
=O(e^{-\epsilon t})
\]
with respect to
$g^{Y,n,0}_{\semiflat|F_{L,\psi}(E,t\theta)}$.
It implies the claim of the theorem.
\hfill\qed

\end{document}

%% file: notation.tex
\newcommand{\nbiga}{\mathcal{A}}

\newcommand{\nbige}{\mathcal{E}}
\newcommand{\nbigf}{\mathcal{F}}

\newcommand{\nbigh}{\mathcal{H}}

\newcommand{\nbigk}{\mathcal{K}}
\newcommand{\nbigl}{\mathcal{L}}
\newcommand{\nbigm}{\mathcal{M}}
\newcommand{\nbign}{\mathcal{N}}
\newcommand{\nbigo}{\mathcal{O}}

\newcommand{\seisuu}{{\mathbb Z}}

\newcommand{\cnum}{{\mathbb C}}
\newcommand{\real}{{\mathbb R}}

\newcommand{\vecA}{{\boldsymbol A}}

\newcommand{\lrarr}{\longrightarrow}

\newcommand{\pf}{{\bf Proof}\hspace{.1in}}
\newcommand{\qed}{\mbox{\rule{1.2mm}{3mm}}}

\def\Hom{\mathop{\rm Hom}\nolimits}

\def\End{\mathop{\rm End}\nolimits}

\def\ad{\mathop{\rm ad}\nolimits}

\def\tr{\mathop{\rm tr}\nolimits}
\def\Tr{\mathop{\rm Tr}\nolimits}

\def\id{\mathop{\rm id}\nolimits}

\newcommand{\del}{\partial}
\newcommand{\delbar}{\overline{\del}}

\def\Def{\mathop{\rm Def}\nolimits}
\def\Harm{\mathop{\rm Harm}\nolimits}

\newcommand{\tildepsi}{\widetilde{\psi}}
\newcommand{\psitilde}{\tildepsi}

\newcommand{\Ltilde}{\widetilde{L}}

\def\hor{\mathop{\rm hor}\nolimits}
\def\ver{\mathop{\rm ver}\nolimits}
\def\Pic{\mathop{\rm Pic}\nolimits}

\def\semiflat{\mathop{\rm sf}\nolimits}
\def\HE{\mathop{\rm HE}\nolimits}
\def\Aut{\mathop{\rm Aut}\nolimits}



%% file: new_theorem.tex
\newtheorem{thm}{Theorem}[section]
\newtheorem{cor}[thm]{Corollary}

\newtheorem{rem}[thm]{Remark}
\newtheorem{lem}[thm]{Lemma}
\newtheorem{prop}[thm]{Proposition}
\newtheorem{df}[thm]{Definition}

%% file: 10-28.bbl
\begin{thebibliography}{99}

\bibitem{Beauville-Narasimhan-Ramanan}
	A. Beauville,
	M. S. Narasimhan,
	S. Ramanan,
	{\em Spectral curves and the generalised theta divisor},
	J. Reine Angew. Math. {\bf 398}, (1989), 169--179.

\bibitem{Dumas-Neitzke}
	D. Dumas,
	A. Neitzke,
	{\em Asymptotics of Hitchin's metric on the Hitchin section},
	Comm. Math. Phys. {\bf 367} (2019), 127--150. 

\bibitem{Fredrickson2}
	L. Fredrickson,
	{\em Exponential decay for the asymptotic geometry
	of the Hitchin metric},
	Comm. Math. Phys. {\bf 375} (2020), 1393--1426. 

\bibitem{Freed}
	 D. Freed,
	 {\em Special K\"{a}hler manifolds},
	 Comm. Math. Phys. {\bf 203} (1999), 31--52.

\bibitem{Gaiotto-Moore-Neitzke}
	D. Gaiotto,
	G.W. Moore,
	A. Neitzke,
	{\em Four-dimensional wall-crossing
	via three-dimensional field theory.}
	Comm. Math. Phys. {\bf 299} (2010), 163--224.

\bibitem{Hitchin-self-duality}
	 N. J. Hitchin,
	 {\em The self-duality equations on a Riemann surface},
	 Proc. London Math. Soc. (3) {\bf 55} (1987), 59--126. 

 \bibitem{MSWW}
	 R. Mazzeo,
	 J. Swoboda,
	 H. Weiss,
	 F. Witt,
	 {\em Ends of the moduli space of Higgs bundles},
	 Duke Math. J. {\bf 165} (2016), 2227--2271. 
 \bibitem{MSWW2}
	 R. Mazzeo,
	 J. Swoboda,
	 H. Weiss,
	 F. Witt,
	 {\em Asymptotic geometry of the Hitchin metric},
	 Comm. Math. Phys. {\bf 367} (2019), 151--191. 
	 
\bibitem{Mochizuki-Asymptotic-Hitchin-metric}
	T. Mochizuki,
	{\em Asymptotic behaviour of the Hitchin metric
	on the moduli space of Higgs bundles},
	arXiv:2305.17638
	
 \bibitem{Mochizuki-improve}
	 T. Mochizuki,
	 {\em Comparison of the Hitchin metric and the semi-flat metric
	 in the rank two case},
	 arXiv:2407.05188.
	 
\bibitem{s1}
C. Simpson,
{\it Constructing variations of Hodge structure
using Yang-Mills theory
and application to uniformization},
J. Amer. Math. Soc. {\bf 1} (1988), 867--918.

\end{thebibliography}
